\documentclass[mn,fleqn]{amsart}
\usepackage{times}
\usepackage{amssymb,amsfonts}
\usepackage[]{graphicx}

\newtheorem{theorem}{Theorem}[section]
\newtheorem{lemma}[theorem]{Lemma}
\newtheorem{corollary}[theorem]{Corollary}

\begin{document}
%%    The information for the title page will be placed between
%%    \begin{document} and \maketitle. The order of most entries
%%    is determined by the class file and can not be changed by
%%    rearranging them. The maketitle command follows after the
%%    abstract.
%%
%%    Most of the following commands will be completed by the publisher.
%%
%%    The copyrightyear is defined in the .clo file as the first argument
%%    of the copyrightinfo command. If the copyrightyear differs from that
%%    value it might be adjusted by the following definition:
%%
%% \renewcommand{\copyrightyear}{2003}% uncomment to change the copyrightyear.
%%
%\DOIsuffix{mana.DOIsuffix}
%%
%% issueinfo for header and copyright line
%\Volume{248}
%\Month{01}
%\Year{2003}
%%
%%    First and last pagenumber of the article. If the option
%%    'autolastpage' is set (default) the second argument may be left empty.
%\pagespan{3}{}
%%
%%    Dates will be filled in by the publisher. The 'reviseddate' and
%%    'dateposted' (Published online) entry may be left empty.
%\Receiveddate{15 November 2003}
%\Reviseddate{30 November 2003}
%\Accepteddate{2 December 2003}
%\Dateposted{3 December 2003}
%%
\keywords{Self homotopy equivalence, aspherical complex, localization.}
\subjclass[msc2000]{55P10, 55P20, 55P60,}

%% \pretitle{Editor's Choice}

%% We have a short and a long form for the title. The short form
%% (optional argument) goes into the running head.

\title[Homotopy equivalences of localized aspherical complexes]{The group of self homotopy equivalences of some
localized aspherical complexes}

%% Please do not enter footnotes or \inst{}-notes into the optional
%% argument of the author command. The optional argument will go into
%% the header.  If there is only one address the marker \inst{x} may be
%% omitted.

%% Information for the first author.
\author[A. Garv\'\i n]{A. Garv\'\i n}\thanks{A. Garv\'\i n, A.
Murillo, and A. Viruel are partially supported by the research
project MTM2004-06262 from the the {\em Ministerio de Educaci\'on
y Ciencia} and by a {\em Junta de Andaluc\'\i a} grant FQM-213}
\address{A. Garv\'\i n, A. Murillo, and A. Viruel\\ Departamento
de \'Algebra, Geometr\'\i a y Topolog\'\i a, Universidad de
M\'alaga\\ AP. 59, 29080 M\'alaga, SPAIN}
%%
%%    Information for the second author
\author[A. Murillo]{A. Murillo}

%%
%%    Information for the third author
\author[J. Remedios]{J. Remedios}\address{J. Remedios\\ Departamento de Matem\'{a}tica Fundamental, Universidad de La Laguna\\
38271 La Laguna, SPAIN.}\thanks{J. Remedios is partially supported
by DGUI of {\em Gobierno de Canarias}.}
%%
%%    Information for the fourth author
\author[A. Viruel]{A. Viruel}
%%
%%    \dedicatory{This is a dedicatory.}
\begin{abstract}

By studying the group of self homotopy equivalences of the localization (at a prime $p$ and/or
zero) of some aspherical complexes, we show that, contrary to the case when the considered space is a
nilpotent complex,  $\mathcal{E}_{\#}^m (X_{p})$ is in general different from $\mathcal{E}_{\#}^m
(X)_{p}$. That is the case even when $X=K(G,1)$ is a finite complex and/or  $G$ satisfies extra
finiteness or nilpotency conditions, for instance, when $G$ is finite or virtually nilpotent.

\end{abstract}
%% maketitle must follow the abstract.
\maketitle                   % Produces the title.

%% If there is not enough space inside the running head
%% for all authors including the title you may provide
%% the leftmark in one of the following three forms:

%% \renewcommand{\leftmark}
%% {First Author: A Short Title}

%% \renewcommand{\leftmark}
%% {First Author and Second Author: A Short Title}

%% \renewcommand{\leftmark}
%% {First Author et al.: A Short Title}

%% \tableofcontents  % Produces the table of contents.
\section{Introduction}

Given a pointed topological space, denote by $\mathcal{E} (X)$, as
usual, the group of homotopy classes of pointed self-homotopy
equivalences of $X$ with the group structure given by composition.
At the present, there is no standard procedure to the
computation of $\mathcal{E} (X)$ and, even for very ``simple'' classes of spaces,
this group is still unknown. Ideas and techniques arising from
completion and localization theory have been applied in different contexts to
approximate $\mathcal{E} (X)$ and some particularly interesting
subgroups.

We shall mainly be concern in one of them, namely the subgroup $\mathcal{E}_{\#}^m (X)$ of
$\mathcal{E} (X)$ formed by those classes inducing the identity  on the homotopy groups up to $m$,
 i.e.,
$$\mathcal{E}_{\#}^m (X)=\ker(\ \mathcal{E} (X)\rightarrow
\Pi_{i=1}^{m}\text{aut}\ \pi_i(X)\ ).$$
Whenever $X$ is either a finite complex, or a space with a finite number of non trivial homotopy
groups, and $m$ is greater than or equal to the ``homotopical" or ``homological" dimension (we denote
it $\dim X$ henceforth), this group is known to be nilpotent \cite{C-M-V,D-Z} as it acts trivially on the
homotopy groups up to $\dim X$. Therefore, it can be localized  in the classical
sense on any set (possibly empty) of
primes $P$ \cite{H-M-R}. What is then
the relation between $\mathcal{E}_{\#}^m (X)_{P}$ and
$\mathcal{E}_{\#}^m (X_{P})$? A complete answer is given whenever
the space is nilpotent \cite{Ma,Pa}. In
fact, for such spaces, the morphism $\mathcal{E}_{\#}^m (X)\rightarrow \mathcal{E}_{\#}^m (X_{P})$, induced
by localization,
 is the
localization morphism and thus $\mathcal{E}_{\#}^m (X_{P})\cong
\mathcal{E}_{\#}^m (X)_{P}$ for any $m\ge\dim X$. In \cite{Mo} the same is proved for completion,
that is to say, $\mathcal{E}_{\#}^m (X_{P}^{\wedge})\cong
\mathcal{E}_{\#}^m (X)_{P}^{\wedge}$.

A natural question arises: can this theorem be extended to a larger class of spaces? the first step in all
the available proofs is proving the result for aspherical complexes. Therefore, we proceed by analyzing the
group
$\mathcal{E}_{\#}^m (X_{p})$ as $X$ runs through  different classes of aspherical spaces, i.e.,
Eilenberg MacLane spaces of the type $K(G,1)$. It is worth to mention that, out of the class of
nilpotent spaces, $p$-localization functors do not, in general,  preserve asphericity, even for virtually
nilpotent $K(G,1)$'s, i.e., when $G$ contains a nilpotent normal subgroup of finite index \cite{B-D,D-M}.
Moreover, for many of these aspherical complexes, $K(G,1)_p$ has an infinite number of non trivial homotopy
groups.

Consider such a space and observe that for any
$m\ge 1$,
$\mathcal{E}_{\#}^m \big(K(G,1)\bigr)=\{1\}$. Hence, realizing a non trivial group as
$\mathcal{E}_{\#}^m \big(K(G,1)_{p}\bigr)$ exhibits a counterexample of the classical results
above for nilpotent spaces.     In this note, we shall do that for aspherical complexes with any extra
nilpotency or finiteness assumptions other than nilpotent of course (finite aspherical complexes and
$K(G,1)$'s with $G$ finite or virtually nilpotent), any prime $p$ and/or zero, and any $m$.

It is known that there are different choices of $p$-localization which extend the classical one for
nilpotent spaces. Except in theorem \ref{kan} below, in which we explicitly indicate that Bousfield
homology localization \cite{B} is being used, we shall be working with the standard $p$-localization with
respect to self maps of $S^1$ \cite{C-P}. In other words, a connected space $X$ is $p$-local if $\pi_1(X)$ is
a
$p$-local group and each $\pi_k(X)$ is a $p$-local $\pi_1(X)$-module, $k\ge 2$.

 The last part of this note heavily rely on basic facts from rational homotopy theory for which we refer to
the standard reference \cite{F-H-T}

\section{Realizing $\mathcal{E}_{\#}^m \big(K(G,1)_{p}\bigr)$}

We shall start by showing that, given a finite nilpotent space $X$, using Bousfield homology
localization
\cite{B} and without any extra assumption on the group $G$, we can easily realize the classical
$p$-localization of the nilpotent group $\mathcal{E}_{\#}^m
(X)$ as a group of the form $\mathcal{E}_{\#}^m \big(K(G,1)_{p}\bigr)$.

\begin{theorem}\label{kan}
Let $X$ be a nilpotent finite complex. Then, there exists a group
$G$ such that, for any $p$ (prime number or zero) and any $m\geq
dimX$,
$$\mathcal{E}_{\#}^m (K(G,1)_{p})\cong \mathcal{E}_{\#}^m
(X)_{p}.$$
\end{theorem}
\begin{proof}
By the Kan-Thurston theorem \cite{K-T} applied to the complex $X$, we find
 a group $G$ and a map $f:K(G,1)\to X$ such that, for any
coefficient system $\mathcal{A}$, the map
$$H^*(f;\mathcal{A}):H^*(K(G,1);\mathcal{A})\to
H^*(X;\mathcal{A})$$ is an isomorphism. In particular, this occurs
when $\mathcal{A}=\mathbb{Z}_{p}$, $\mathbb{Q}$, and therefore,
the Bousfield homology $p$-localization ($p$ a prime number or
zero) of $K(G,1)$ coincides with that of $X$, $K(G,1)_{p}\simeq
X_{p}$.

However, the Bousfield localization is an extension of
 the classical localization for nilpotent spaces. Therefore, the theorem of Maruyama stated above
\cite{Ma,Pa} implies that, for any $m\geq dim X$, $\mathcal{E}_{\#}^m (X_{p})\cong \mathcal{E}_{\#}^m
(X)_{p}$ and therefore $\mathcal{E}_{\#}^m (K(G,1)_{p})\cong
\mathcal{E}_{\#}^m (X)_{p}$.
\end{proof}

\begin{corollary}
Let $H$ be any finite abelian $p$-local group without $2$-torsion. Then, for any $m>2$,
there exists a group $G$  such that $$\mathcal{E}_{\#}^m(K(G,1)_{p})\cong H.$$
\end{corollary}

\begin{proof}
Indeed, by \cite[Remark 4.7 (4)]{A-M},
 any finite abelian group without
$2$-torsion $H$ is isomorphic to $\mathcal{E}_{\#}^m(X)$ in which $X$ is a
co-Moore space of dimension $m>2$ and of type $(\mathbb Z\oplus H, m-1)$. To finish apply theorem above
to $X$ taking into account that $H_{p}=H$.
\end{proof}

This is our first class of spaces of the form $K(G,1)$ for which
$\mathcal{E}_{\#}^m(K(G,1)_{p})$ is not the trivial group for
some $p$ and some $m\ge 1$. However, the group $G$ obtained from Kan-Thurston theorem is far
from satisfying finiteness or nilpotency restrictions. In the next
class of examples we impose this kind of restricted
behavior. We begin by the following observation:

\begin{lemma}\label{lema}
Let $G$ be a perfect, virtually nilpotent group. Then for any
prime $p$
$$\mathcal{E}_{\#}^1(K(G,1)_{p})=
\mathcal{E}(K(G,1)_{p}).$$
\end{lemma}
\begin{proof} Recall \cite{B-C} that a group $G$ is called generically trivial if, for
any prime $p$, $G_p=\{1\}$. It has been proved, first for finite
groups \cite{B-C} and then for virtually nilpotent groups \cite{D-M}, that generically
trivial groups coincide with perfect groups. Hence, taking into
account that $\pi_1(K(G,1)_{p})=G_p=\{1\}$ \cite{C-P}, any self-equivalence
of $\mathcal{E}(K(G,1)_{p})$ is trivially in
$\mathcal{E}_{\#}^1(K(G,1))_{p}$.
\end{proof}

\begin{theorem} For any group $H$ and any prime $p$, there exists a group
 $G$ such that $H$ is a subgroup of
$\mathcal{E}_{\#}^1(K(G,1)_{p})$. Moreover, if $H$ is a finite group, then $G$ can also be chosen finite.
\end{theorem}
\begin{proof} Consider the ``regular" representation of $H$ as a subgroup of the group $\Sigma(H)$ of
bijections of
$H$. Choose any  group $F$ for which the space $K(F,1)_{p}$ is non
contractible and simply connected (i.e. $F_p=\{1\}$). For instance, any group as in lemma above satisfies this for all
$p$. In
\cite{B-D,D-M}, the reader may find many examples of these groups which satisfy additional finiteness or nilpotency
properties. Under these assumptions consider the free product $G=*_{h\in H}F$ and observe that $BG=K(G,1)=\vee_{h\in
H}K(F,1)=\vee_{h\in H}BF$. Hence
$K(G,1)_p=\bigl(\vee_{h\in H}\,K(F,1)\bigr)_p$ and therefore, via  Van Kampen theorem, this
space is also 1-connected and  non contractible as, for instance, it has   non trivial mod $p$ homology. Next,
consider the composition
$$
\gamma\colon
 \Sigma(H){\stackrel{\varphi}{\hookrightarrow}}\mathcal{E}\bigl(\vee_{h\in
H}\,K(F,1)_p\bigr){\stackrel{\psi}{\longrightarrow}}\mathcal{E}\bigl(K(G,1)_p\bigr)
$$
where $\varphi$ is the injective morphism defined by $\varphi(\alpha)(x_h)_{h\in H}=(x_{\alpha(h)})_{h\in H}$ and $\psi$
is the map $\mathcal{E}(X)\to\mathcal{E}(X_p)$ induced by localization. We now show that $\gamma$ is an injective
morphism: indeed, for each $h\in H$ consider the canonical inclusion $i_h\colon F\to G$ and ``projection" $q_h\colon
G\to F$ so that $q_h\circ i_h=1_F$ and $q_h\circ i_{h'}$ is the trivial map for $h\not=h'$. Now, assume there is a
bijection $\alpha\in\Sigma(H)$, $\alpha\not=1_H$, such that $\gamma(\alpha)=1_{K(G,1)_p}$. Hence, there exists $h\in H$
for which $\alpha(h)=h'\not=h$. For this element we have:
$$
1_{K(F,1)_p}=(B1_F)_p=(Bq_h)_p\circ (Bi_h)_p=(Bq_h)_p\circ\gamma(\alpha)\circ (Bi_h)_p=(Bq_h)_p\circ(Bi_{h'})_p=*
$$
However, $K(F,1)_p$ is non contractible and therefore $\gamma$ is an injective morphism. To finish, consider the
restriction of $\gamma$ to $H$,
$$
H\hookrightarrow K(G,1)_{p},
$$
 an observe that, since $K(G,1)_p$ is
simply connected,
$\mathcal{E}_{\#}^1(K(G,1)_{p}) =\mathcal{E}(K(G,1)_{p})$.

On the other hand, note that the group $G$ is not finite even when $H$ is. Hence, whenever this is the case, to prove the
second part of the theorem, we follow the same argument choosing the finite group
$G=\prod_{h\in H} F$. Again, $K(G,1)_p=\prod_{h\in H}\,K(F,1)_p$  is also 1-connected
and  non contractible.
As before, consider the restriction to $H$ of the injective morphism
$$
 \Sigma(H)\hookrightarrow \mathcal{E}(K(G,1)_{p})\cong
\mathcal{E}(\prod_{h\in H}K(F,1)_{p})
$$
\end{proof}

The aspherical spaces studied up to this point are not, in general, finite complexes. To cover this case,
we now present a
$K(G,1)$ which is a  virtually nilpotent (i.e, $G$ is virtually nilpotent) finite complex
whose rationalization has non trivial self homotopy equivalences which fix the homotopy groups up to
any dimension. Recall that, from the homotopical point of view, an { \em infranilmanifold} is a manifold of
the homotopy type of an aspherical complex $K(G,1)$ in which $G$ is a torsion free, finitely generated,
virtually nilpotent group containing no non trivial finite normal subgroups \cite{De,Ch}. In fact, many of
such spaces are flat Riemannian manifolds.

\begin{theorem}\label{teorema}
There exist a compact infranilmanifold of the form $K(G,1)$ for which
$\mathcal{E}_{\#}^m(K(G,1)_{\mathbb{Q}})\neq\{1\}$ for all $m\geq dim\
K(G,1)$.
\end{theorem}
\begin{proof} Let $K(F,1)$ be the non orientable $4-$manifold \cite{Ch,B-D} in which $F$ is the virtually nilpotent group
generated by $\{x_1,x_2,x_3,x_4,\alpha,\beta\}$ with relations:
\begin{itemize}

\item $[x_i,x_j]=1,\ i\neq j,$

\item $\alpha^2=x_3,\ \beta^2=x_4,\
\alpha\beta=x_2^{-1}x_3x_4^{-1}\beta\alpha,$

\item $\alpha x_1=x_1^{-1}\alpha,\ \alpha x_2=x_2^{-1}\alpha,\ \alpha
x_3=x_3\alpha,\ \alpha x_4=x_4^{-1}\alpha,$

\item $\beta x_1=x_1^{-1}\beta,\ \beta x_2=x_2^{-1}\beta,\ \beta
x_3=x_3^{-1}\beta,\ \beta x_4=x_4\beta.$
\end{itemize}
It is indeed  an infranilmanifold of dimension $4$ whose
rationalization turns out to be \cite[Example 5.2]{B-D}:
$$K(F,1)_{\mathbb{Q}}\simeq (S^2\vee
S^3\vee S^3)_{\mathbb{Q}}.
$$
Define $G=F\times \mathbb{Z}$ and observe that
$K(G,1)=K(F,1)\times S^1$ is in fact a $5-$dimensional manifold
whose rationalization is
$$K(G,1)_{\mathbb{Q}}=(S^2\vee
S^3\vee S^3)_{\mathbb{Q}}\times S_{\mathbb{Q}}^1.$$
We shall prove that
$\mathcal{E}_{\#}^m(K(G,1))\neq
\{1\}$ for all $m$. For that observe in the first place that
$(S^2\vee S^3\vee S^3)\times S^1$ is a nilpotent space so that, in
view of \cite{Ma},
$$\mathcal{E}_{\#}^m\bigl((S^2\vee
S^3\vee S^3)_{\mathbb{Q}}\times S_{\mathbb{Q}}^1\bigr)\cong \mathcal{E}_{\#}^m(( S^2\vee S^3\vee
S^3)\times S^1)_{\mathbb{Q}}.
$$
Moreover, this group coincides with
the group $\mathcal{E}_{\#}^m(\Lambda V,d)$ defined as
follows:
\medskip

 (1) $(\Lambda V,d)$ is the minimal model of the space
$( S^2\vee S^3\vee S^3)\times S^1$. In other words, $\Lambda V$ is the commutative
free algebra generated by the graded vector space $V\cong \pi_*((
S^2\vee S^3\vee S^3)\times S^1)\otimes \mathbb{Q}$, and $d$ is
a certain differential satisfying $dV\subset \Lambda ^{\geq 2}V$, i.e.,
for each generator $v\in V$, $dv$ is a polynomial in $\Lambda V$
with no linear terms.
\medskip

 (2) $\mathcal{E}_{\#}^m(\Lambda V,d)$ is the subgroup of homotopy
classes of differential graded algebra automorphisms  $f\colon(\Lambda
V,d)\to (\Lambda V,d)$ which satisfy $f(v)-v\in\Lambda^{\geq 2}V$ for
any $v\in V$ of degree less than or equal to $m$.
\medskip

To compute $\mathcal{E}_{\#}^m(\Lambda V,d)$ observe in the first place that
$(\Lambda V,d)=(\Lambda W,d)\otimes (\Lambda x,0)$ where
$(\Lambda W,d)$ is the minimal model of $( S^2\vee S^3\vee S^3)$
and $(\Lambda x,0)$, where $x$ is a single element of degree $1$,
is the minimal model of $S^1$. On the other hand, it is well known \cite{F-H-T} that the space $(
S^2\vee S^3\vee S^3)$ is both formal and coformal, that is to say, its rational homotopy type
depends only on either, its rational cohomology algebra, or its rational homotopy Lie algebra. To fix
notation  let
$H=H^*( S^2\vee S^3\vee S^3;\mathbb{Q})$ be the commutative algebra
generated by an element $a_2$ of degree $2$ and two elements $b_3, c_3$ of degree $3$, with trivial
multiplication. Then, the model $(\Lambda W,d)$ of $( S^2\vee S^3\vee S^3)$ is classically
constructed as follows
\cite{F-H-T}: the vector space $W=W_*^*$ is built inductively bigraded and this  bigrading is inherited
by $\Lambda W$, according to the usual rule,
 so that it satisfies:
\begin{itemize}

\item $W_0\cong H$, $dW_0=0$,
\item $W_{m+1}\cong (\Lambda^2W_{\leq m})_m\cap ker\ d$;\ \
$d:W_{m+1}\stackrel{\cong}{\longrightarrow}(\Lambda^2W_{\leq
m})_m$,
\item $H^*(\Lambda W,d)\cong H$,
\item Any cohomology class in $H^*(\Lambda W,d)$ represented by a
decomposable cycle vanishes.
\end{itemize}

Our theorem will then be established once we prove the following:

\begin{lemma}

$\mathcal{E}_{\#}^m(\Lambda V,d)\neq \{1\}$ for all $m$.
\end{lemma}

\medskip\noindent We shall define an automorphism $\varphi$ of $(\Lambda
V,d)=(\Lambda W,d)\otimes (\Lambda x,0)$, inductively on $W_m$,
satisfying:
\begin{itemize}

\item $\varphi (a_2)=a_2,$ $\varphi (b_3)=b_3,$ $\varphi
(c_3)=c_3+xa_2,$ $\varphi (x)=x$,
\item For each $w$ generator of $W_m$, $m\geq 0$, $\varphi
(w)=w+w'x$, $w'\in W$.
\end{itemize}

Obviously $\varphi _{|_{W_0}}$ satisfy our hypothesis.  Assume
it has been extended to $\Lambda W_{\leq m-1}\otimes\Lambda x$ and let $w\in W_m$. Then,
$dw\in(\Lambda^2W_{\leq m-1})_{m-1}$, i.e.,
$\displaystyle dw=\sum_i u_iv_i$, $u_i,v_i\in W_{\leq m-1}$.
Hence, $\displaystyle\varphi (dw)=\sum_i\varphi (u_i)\varphi
(v_i)=\sum_i(u_i+w_i'x)(v_i+w_i''x)= \sum_iu_iv_i+
\sum_i(u_iw_i''\pm w_i'v_i)x= dw+\alpha x$, being $\displaystyle
\alpha= \sum_i(u_iw_i''\pm w_i'v_i)$. On the other hand $\varphi
(dw)$ is obviously a cycle so is $\alpha$. But, since $\alpha$ is
decomposable, it must be a boundary: $\alpha=dw'$,
$w'\in W$. Thus
$$
\varphi (dw)= d(w+w'x).
$$
Finally define $\varphi (w)=w+w'x$ and observe that
the automorphism $\varphi $ of $(\Lambda V,d)$ just defined is not homotopic to the
identity as $H^*(\varphi )\neq 1_{H^*(\Lambda V,d)}$. Indeed
$H^*(\varphi )[c_3]=[c_3]+[x][a_2]$. This finishes the proof of the lemma and thus, that of
theorem \ref{teorema}.
\end{proof}

\nocite{*}

\bibliographystyle{plain}

\end{document}